\documentclass[11pt,noamsfonts]{amsart}

\usepackage{amsmath,amsfonts,amssymb,amscd,verbatim,xypic, xy}

\newcommand{\marginlabel}[1]%
  {\mbox{}\marginpar{\raggedleft\hspace{0pt}\bfseries\sf#1}}

\def\ZZ{{\mathbb Z}}

\def\QQ{{\mathbb Q}}
\def\PP{{\mathbb P}}

\def\cI{\mathcal{I}}
\def\cJ{\mathcal{J}}
\def\cD{\mathcal{D}}
\def\cA{\mathcal{A}}

\def\cL{\mathcal{L}}

\def\cN{\mathcal{N}}

\DeclareMathOperator{\Ker}{ker}

\newtheorem{lemma}{Lemma}[section]
\newtheorem{theorem}[lemma]{Theorem}
\newtheorem{corollary}[lemma]{Corollary}
\newtheorem{proposition}[lemma]{Proposition}
\newtheorem*{claim}{Claim}

\theoremstyle{definition}
\newtheorem{definition}[lemma]{Definition}
\newtheorem{example}[lemma]{Example}

\newtheorem{remark}[lemma]{Remark}

\newtheorem*{notation}{Notation}

\numberwithin{equation}{section}

\newcommand{\bean}{\begin{eqnarray}}
\newcommand{\eean}{\end{eqnarray}}
\newcommand{\be}{\begin{displaymath}}
\newcommand{\ee}{\end{displaymath}}
\newcommand{\bea}{\begin{eqnarray*}}
\newcommand{\eea}{\end{eqnarray*}}

\newcommand{\ol}{\overline}

\xyoption{arc}

\begin{document}

\title{Tautological rings of stable map spaces}

\author[Anca M. Musta\c{t}\v{a}]{Anca~M.~Musta\c{t}\v{a}}
\author[Andrei Musta\c{t}\v{a}]{Andrei~Musta\c{t}\v{a}}
\address{Department of Mathematics, University of Illinois, 1409 W. Green Street, Urbana, IL 61801, USA} \email{{\tt amustata@math.uiuc.edu,
dmustata@math.uiuc.edu}}

\date{\today}

\begin{abstract}
 We find a set of generators and relations for the system of extended 
tautological rings associated to the  moduli spaces of stable maps in genus 
zero, admitting a simple geometrical interpretation. In particular, when the 
target is $\PP^n$, these give a complete presentation for the cohomology and 
Chow rings in the cases with/without marked points. 
\end{abstract}
\maketitle

\bigskip

\section*{Introduction}

 Let $d\in H_2(X)$ be a curve class on a smooth projective variety $X$. The space 
$\ol{M}_{0,0}(X,d)$ parametrizes maps from rational smooth or nodal curves into $X$ with image 
class $d$, such that any contracted component contains at least 3 nodes. Over $\ol{M}_{0,0}(X,d)$ 
there exists a tower of moduli spaces of stable maps with marked points and morphisms
$\ol{M}_{0,m+1}(X,d)\to \ol{M}_{0,m}(X,d)$ forgetting one marked point, such that $\ol{M}_{0,m+1}(X,d)$, together with an evaluation map $ev_{m+1}: \ol{M}_{0,m+1}(X,d)\to X$ and $m$ natural sections $\sigma_i$, 
 form the universal family over $\ol{M}_{0,m}(X,d)$. 

 In this paper we investigate the relation between the structure of the cohomology ring of the variety
 $X$ and that of $\ol{M}_{0,0}(X,d)$, which is less than obvious in particular when $m=0$. When marked points
 exist, pullback by the natural evaluation maps $ev_i:\ol{M}_{0,m}(X,d)\to X$ generate a set of classes 
on the moduli space. A system of tautological rings for $\{\ol{M}_{0,m}(X,d)\}_m$ is constructed from 
these classes by analogy with the moduli space of curves. When the target space is an $SL_n$--flag varity, these coincide with the 
cohomology rings as shown in (\cite{oprea1}). We expect that the same result holds for a larger class of varieties. The Gromov-Witten invariants of $X$ and 
their gravitational descendants are intersection numbers in these rings. Special tautological classes
are the first Chern classes $\psi_i$ of the tautological line bundles over $\ol{M}_{0,m}(X,d)$. Another
 special set of tautological classes are the $k$-- classes $k_a(\alpha)$ defined in \cite{kk} by the formula
$k_a(\alpha)=f_*(\psi_i^{a+1} ev_i^*\alpha )$, where $f$ is the forgetful map above.

We define a set of extended tautological rings, offering
a more convenient encoding of the boundary, and for which the original tautological rings are invariant subrings. Theorems 2.5  and 2.6 in 
the text identify a set of generators and relations for these rings in terms of boundary strata and $k$-- classes
$k_{-1}(\alpha_i)$ of the generators $\{\alpha_i\}_i$ for the ring $H^*(X)$. Relations in $H^*(X)$ 
induce relations in the tautological rings via $k$-- class decomposition. Natural universal 
relations exist on the boundary. 

In particular, when $X=\PP^n$, the above form a complete presentation of the extended cohomology 
(and Chow) rings of $\ol{M}_{0,0}(\PP^n,d)$, with/ without marked points, of a simple geometrical
interpretation. (Theorem 3.3 in text).  Thus we recover the degree 3 and 2 cases formulated by Behrend and O'Halloran in \cite{behrend2}.

\section{Definition of Tautological rings}

For any smooth projective target $X$, the moduli spaces $\ol{M}_{0,m}(X,d)$ parametrising pointed stable maps of class $d$ into $X$ are Deligne-Mumford stacks sharing a couple of common features with the moduli spaces of stable curves:
\begin{enumerate}
\item natural forgetful morphisms $f: \ol{M}_{0,m+1}(X,d) \to \ol{M}_{0,m}(X,d)$;
\item natural gluing maps $\ol{M}_{\tau}\to  \ol{M}_{0,m}(X,d)$.
  \end{enumerate}
Here $\ol{M}_{\tau}$ is a fiber product of moduli spaces $ \ol{M}_{0,m_i}(X,d_i)$ along evaluation morphisms to $X$, such that $\sum_id_i=d$ and $\ol{M}_{\tau}$ represents maps into $X$ from stable reducible curves of a split type $\tau$. We say that $\ol{M}_{\tau}$ is a boundary stratum.
 
 The unique features of  $\ol{M}_{0,m}(X,d)$ consist in the existence of $m$ evaluation maps into $X$, and the existence of virtual fundamental classes $[\ol{M}_{0,m}(X,d)]^{vir}$ in  Chow groups and homology (\cite{behrend-fantechi}). These classes are compatible with the natural morphisms above (\cite{behrend-manin}) and compensate for the fact that  $\ol{M}_{0,m}(X,d)$ in general is not a smooth stack of the expected dimension.

With these we can construct  a system of rings inside $H_*(\ol{M}_{0,m}(X,d))$. These are called the tautological rings, and are defined by analogy with the moduli space of curves.

We will work in homology of $\ol{M}_{0,m}(X,d)$ via  the homomorphism
$$\sigma: H^*(\ol{M}_{0,m}(X,d))\to H_*(\ol{M}_{0,m}(X,d))$$
$$\sigma(\beta ) = \beta\cap [\ol{M}_{0,m}(X,d)]^{vir}.$$

The classes $\sigma( ev_i^*\alpha)$ and their pushforwards via forgetful and gluing maps are the first examples of tautological classes. 
 We note that the forgetful maps  admit natural flat
pullback in homology.  The gluing maps $j$ of the boundary strata admit a Gysin map, which will be denoted by $j^*$.

 Intersection of tautological classes in $H_*(\ol{M}_{0,m}(X,d))$ can be defined in a canonical way. Consider three moduli spaces $\ol{M}_1$, $\ol{M}_2$, $\ol{M}$ related by natural morphisms $f_1$, $f_2$, and the fibre square diagram
\bea  \xymatrix{  \ol{M}_1\times_{\ol{M}}\ol{M}_2 \ar[d]^{\tilde{f}_2}\ar[r]^{\tilde{f}_1}    &  \ol{M}_2 \ar[d]^{f_2}\\
 \ol{M}_1 \ar[r]^{f_1} & \ol{M}    }  \eea
and let $f:=f_1\circ \tilde{f}_2=f_2\circ\tilde{f}_1$. 
Here $f_1$ and $f_2$ are strata embeddings, forgetful maps or compositions of these, or the identity. Similarly we define $$ \sigma : H^* ( \ol{M}_i ) \to H_*( \ol{M}_i ) $$  $$\sigma(\beta ) = \beta\cap [\ol{M}_i]^{vir},$$
for $i=1,2$.

\begin{definition}
The intersection product of two classes $f_{1 *}(\sigma (\beta_1))$ and $f_{2 *}(\sigma (\beta_2))$ is defined as
$$f_{1 *}(\sigma (\beta_1))\cdot f_{2 *}(\sigma (\beta_2)):=f_*(\tilde{f}_1^*\beta_2\cdot\tilde{f}_2^*\beta_1\cap \tilde{f}_1^*([\ol{M}_2]^{vir})) $$
\end{definition}
 Note that by the axioms of the virtual class, $\tilde{f}_1^*([\ol{M}_2]^{vir})= \tilde{f}_2^*([\ol{M}_1]^{vir})=f^*([\ol{M}]^{vir})$, so the product is well defined. The product is clearly associative. We may thus define

\begin{definition}
 The tautological system of rings of $X$, $d$ is the set of smallest $\QQ$-- algebras inside the cohomology groups
 $$R(\ol{M}_{0,m}(X,d)) \subset H_*(\ol{M}_{0,m}(X,d))$$
satisfying the following properties:
\begin{enumerate}
\item $R(\ol{M}_{0,m}(X,d))$ contains all classes $\sigma(ev_i^*(\alpha))$, where $\alpha\in H^*(X)$, $i=1,...,m$.
\item The system is closed under push-forward via all forgetting maps: $$ f_*:R(\ol{M}_{0,m+1}(X,d))\to R(\ol{M}_{0,m}(X,d)).$$
\item The system is closed under push-forward via all gluing maps:
$$j_*:  R(\ol{M}_{0,A_1\cup \{\bullet\}}(X,d_1))\otimes_{\QQ} R(\ol{M}_{0,\{\bullet\}\cup A_2  }(X,d_2))\to R(\ol{M}_{0,m}(X,d))$$
where $A_1\bigcup A_2=\{1,...,m\}$ and $d_1+d_2=d$.
\item The system is closed under the product defined above.
\end{enumerate}
\end{definition}
\subsection{}

Consider a smooth projective variety $X$ and an embedding $X\hookrightarrow \prod_{i=1}^s\PP^{n_i}$ such that the algebraic part of $H^2(X,\ZZ)$ is $\ZZ^s$, generated by first Chern classes of the very ample line bundles $\cL_1, ...,\cL_s$ giving the embedding above. Thus every curve class $d$ on $X$ is Poincar\'e dual to the first Chern class of a line bundle $\otimes_{i=1}^s \cL_i^{d_i}$ for nonnegative integers $d_i$. We write $d=(d_1,...,d_s)\in \ZZ^s$. Let $G=S_{d_1}\times...\times S_{d_s}$ denote the product of $s$ groups of permutations $S_{d_i}$.  We recall succinctly the construction of a $G$-network of gluing morphisms for $\ol{M}:=\ol{M}_{0,m}(X,d)$, and its extended homology groups (see \cite{noi1} and \cite{noi2} for a more detailed presentation).

 Let $I$ be a set whose elements are of the form $h=h_1\sqcup ...\sqcup h_s\sqcup M_h$ such that each $h_i\subset \{1,...,d_i\}$ and $M_h\subset\{2,...,m\}$.
 Assume that  $h\cap h'=h, $ $h'$ or $\emptyset$  for all $ h, h'\in I$. 
\begin{definition}
The space $\ol{M}_I$
parametrizes degree $d$  stable maps $\varphi: (C, \{p_j\}_{j=1,...,m})\to X $
together with a set of closed curves $ \{C_{h}\}_{h\in I}$ such that:
\begin{enumerate}
\item $\forall h\in I$, $p_1\not\in C_h\subset C$, and the degree of the map $\varphi_{| C_{h}}$ is $(|h_1|,...,|h_s|)$.
\item the incidence relations among the elements of $I$ translate into analogous incidence relations among the curves $C_h$:
\begin{itemize}
\item $\forall h\in I$, $\forall i\in \{2,...,m\}$,  $p_i\in C_h$ iff $i \in M_h$.
\item $C_h\subset C_{h'}$ iff $h\subset h'$ and $C_h\bigcap C_{h'}=\emptyset$ iff $h\bigcap h'=\emptyset$.
\end{itemize}
 \end{enumerate}
 \end{definition}

 In particular, $\ol{M}=\ol{M}_{\emptyset}$. For any two subsets $I\subset J$, there is natural local regular embedding $\phi^I_J: \ol{M}_{J}\to \ol{M}_{I}$.  The spaces $\ol{M}_{I}$ with the morphisms $\phi^I_J$ form a network, on which the group $G$ acts naturally. 

 In the case of $\ol{M}= \ol{M}_{0,0}(X, d)$, let  $I$ denote any set of 2--partitions $h\bigcup \bar{h}= \bigsqcup_i \{1,...,d_i\}$ such that for any pair $(h,\bar{h}),$ $(h',\bar{h}') \in I$, the set $\{h\bigcap h', h\bigcap\bar{h}', \bar{h}\bigcap h', \bar{h}\bigcap \bar{h}'\}$ has exactly three elements. $\ol{M}_{I}$ represents stable maps $\varphi: C\to X$, together with marked splittings $C_{h}\bigcup C_{\bar{h}}=  C$ for all $(h,\bar{h})\in I$ satisfying natural incidence relations as in Definition 1.2.

 For each set $I$ as above, let $G_I\subset G$ be the subgroup which fixes all elements of $I$. Any
$g\in G$ induces a canonical isomorphism
$g:\ol{M}_I\to\ol{M}_{g(I)}$.

\begin{definition} The extended homology 
 groups are \bea
B_l(\ol{M}_{0,m}(X,d))  :=\oplus_{I} H_{l}(\ol{M}_{I}) / \sim , \eea  the
sum taken after all subsets $I$ as above. The equivalence relation $\sim
$ is generated by: $$\phi^I_{J*}(\alpha)\sim \sum_{[g]\in
G_I/G_J} g_*(\alpha)$$ for any $\alpha \in H_*(\ol{M}_{J})$ and
$J\supset I$.


\end{definition}

A system of extended tautological rings is defined inside $B_*(\ol{M}_{0,m}(X,d))$:
\begin{definition}
The system of extended tautological rings of $X$, $d$ is the set of smallest $\QQ$-- algebras inside the cohomology groups
 $$T(\ol{M}_{0,m}(X,d)) \subset B_*(\ol{M}_{0,m}(X,d))$$
satisfying all the properties listed in Definition 1.1, where the  maps $\phi^I_J$ stand for gluing morphisms.

\end{definition}

  In \cite{noi1}, \cite{noi3} we give a more detailed account of the extended intersection theory associated to the  networks of boundary maps in a moduli space. $T(\ol{M}_{0,m}(X,d))$ admits a natural multiplicative structure, and enjoys the salient properties of $R(\ol{M}_{0,m}(X,d))$, like existence of pushforward and pullback for the forgetful maps and for the gluing morphisms. The most important feature of  $T(\ol{M}_{0,m}(X,d))$ is that it allows the fundamental classes of boundary strata to be decomposed as polynomials of boundary divisors. This greatly simplifies the task of listing ring generators. On the other hand,  $R(\ol{M}_{0,m}(X,d))$ is the invariant subalgebra of  $T(\ol{M}_{0,m}(X,d))$  with respect to the natural action of $G$. In \cite{noi2} we show how  an additive basis for  $R(\ol{M}_{0,m}(\PP^n,d))$  may be extracted from the ring structure of  $T(\ol{M}_{0,m}(\PP^n,d))$ in the case $m>0$.

When $X$ is convex, the moduli spaces are smooth and thus the homomorphism $\sigma$ is an isomorphism. Then $R(\ol{M}_{0,m}(X,d))$ may be viewed naturally as a $\QQ$--subalgebra of the bivariant cohomology $H^*(\ol{M}_{0,m}(X,d))$ and its definition coincides with the one introduced in \cite{oprea1}.

 In the case of $SL_N$-flag varieties, it was shown in \cite{oprea1} that the cohomology rings are tautological, and isomorphic to the Chow rings. In general we may ask for which set of varieties the same property is true. With the techniques of \cite{noi1} and \cite{noi3}, whether the cohomology of stable map spaces is tautological depends essentially on the structure of the linear sigma model (quasi-map space) of $X$. We may expect the answer to be affirmative in other cases like smooth toric varieties.

\section{Structure of tautological rings}

The set of stable map spaces for a target $X$ may be refined by considering a larger variety of stability conditions.

\begin{notation}
Let $A$, $B$ be two disjoint sets such that $A\bigcup B=\{1,...,m+1\}$ and let $d_A, d_B, d$ be curve classes on $X$ such that $d_A+d_B=d$. We denote by $D(A,d_A|B,d_B)$, the divisor of $\ol{M}_{0,m}(X,d)$ representing split curves, with a degree $d_A$-- component containing the set $A$ of marked points, and a degree $d_B$-- component containing the set $B$ of marked points. Assume $1\in B$. With the notations from Definition 1.3, the support of  $D(A,d_A|B,d_B)$ is the image of the morphism $\ol{M}_h\to \ol{M}_{0,m}(X,d)$ if $m>1$, or  $\ol{M}_{h\bar{h}} \to \ol{M}_{0,0}(X,d)$, where $h=h_1\sqcup ...\sqcup h_s\sqcup A$ is such that $(|h_1|,...,|h_s|)=d_A$. 
\end{notation}

  In \cite{noi3} we constructed a sequence of Deligne-Mumford stacks $\ol{U}^k_{m}$, and of morphisms
$$\ol{M}_{0,m+1}(X,d)=\ol{U}^0_m\to...\to \ol{U}^{l_m}_m\to \ol{M}_{0,m}(X,d),$$ such that each $\ol{U}^k_m$ is birational to $\ol{M}_{0,m+1}(X,d)$. When $m=0$ we drop the subscript $m$.
  Denote the morphisms above and their compositions by $f^{k'}_k: \ol{U}^{k'}_m\to \ol{U}^k_m$ for any $ k'\leq k$,  $f^{k}: \ol{U}^{k}_m\to \ol{M}_{0,m}(X,d)$. Each step $f^{k+1}_k$ is a blow-up along a regularly embedded codimension two locus, such that pullback  to $\ol{M}_{0,m+1}(X,d)$ of the exceptional divisor is one of the divisors $D(A,d_A|B,d_B)$ above, with  $1\in B$ (if $m>0$) and $m+1\in A$.
In the case $m=0$, for every pair $(d_A, d_B)$ as above, exactly one of the divisors $D(\{ 1\},d_A|\emptyset, d_B)$,  $D(\emptyset, d_A|\{ 1\}, d_B)$ on $\ol{M}_{0,1}(X,d)$ is an exceptional divisor in the sequence of contractions above, except in the case $d_A=d_B$. The choice of the exceptional divisor within the pair depends on the nature of the stability conditions imposed at each step. 


 In general,  $f^{l_m}:  \ol{U}^{l_m}_m\to \ol{M}_{0,m}(X,d)$ is a $\PP^1$-- bundle. The only  exception is when $m=0$ and there exists a class $d'$ such that $2d'=d$. Then the map $f^{l_m}$ may be formally understood at the level of \'etale covers as the a blow-up of a $\PP^1$-- bundle.
( see \cite{noi3} for a complete argument). 



 The morphisms making up the factorization of the forgetful map are locally complete intersection morphisms, so they admit pullback in homology. Consequently, tautological and extended tautological rings $T( \ol{U}^k_m)$ may be defined in a natural way for each of the contractions $\ol{U}^k_m$, via $[\ol{U}^k_m]^{vir}:=f^{k *}[\ol{M}_{0,m}(X,d)]^{vir}$. 

Let $m=0$. To fix notation for strata in the extended ring, consider for each step $k$ all the sets $h$ as in Definition 1.3, such that the following diagram consists of two fiber squares
 \bea \xymatrix{  \ol{U}^{k-1}_{h}\bigcup\ol{U}^{k-1}_{\bar{h}} \ar[d]_{j_h^{k-1}, j_{\bar{h}}^{k-1}} \ar[r]^{  f^{k-1}_{k h\bar{h}} } & \ol{U}^{k}_{h\bar{h}}  \ar[d]^{j^k_{h\bar{h}}} \ar[r]^{ f^{k}_{h\bar{h}}} &  \ol{M}_{h\bar{h}} \ar[d]^{j_{h\bar{h}}} \\
\ol{U}^{k-1} \ar[r]^{ f^{k-1}_k} & \ol{U}^k  \ar[r]^{ f^{k}} & \ol{M}_{0,0}(X,d).
 }\eea
 Here  vertical morphisms are regular local embeddings of boundary strata. The generic fiber of $f^{k}_{h\bar{h}}$ is an irreducible curve, while  $f^{k-1}_{k h\bar{h}}$ consists of two morphisms 
$f^{k-1}_{k h}$ and $f^{k-1}_{k \bar{h}}$, each with irreducible generic fiber.  
$f_{h\bar{h}}^{k+1}$  admits a section $s_h:
\ol{M}_{h\bar{h}}\to \ol{U}^k_{h\bar{h}}$. Let   $S_h$ denote the
image of $s_h$, with regular embedding  $j: S_h\to \ol{U}^{k}_{h\bar{h}}$, and $j^{k}_{h}:= j_{h\bar{h}}^{k}\circ j$. The  space $\ol{U}^{k-1}$ is  the blow-up of $\ol{U}^{k}$ along
$S_h$, and $\ol{U}^{k-1}_h$ is its exceptional
divisor, with  projection $g_h: \ol{U}^{k-1}_h\to  S_h$. Let $f^{k-1}_h= f^{k}_{h}\circ f^{k-1}_{k h}$. The space $\ol{U}^{k-1}_{\bar{h}}$ is the strict transform of
$\ol{U}^{k}_{h\bar{h}}$  and $f^{k-1}_{\bar{h}}:= f^{k}_{h\bar{h}} \circ f^{k-1}_{k \bar{h}}$. The two divisors $\ol{U}^{k-1}_h$ and $\ol{U}^{k-1}_{\bar{h}}$ intersect along $\tilde{S}_h\cong S_h$.

When $k=1$ i.e. $ \ol{U}^{k-1} = \ol{M}_{0,1} (X,d)$ we drop the superscript $k-1$, writing simply $f_h$, $f_{ \ol{h}}$, $j_h$, $j_{\ol{h}}$ for the appropriate morphisms.

  The fundamental class of $\ol{U}^{k}_{h\bar{h}}$ in  $T(\ol{U}^{k})$ will be denoted by $D_{h\bar{h}}$, and is the pullback of the analogous class on $\ol{M}_{0,0}(X,d)$. The classes of $\ol{U}^{k-1}_h$ and $\ol{U}^{k-1}_{\bar{h}}$ in $\ol{U}^{k}$ are $D_h$ and $D_{\bar{h}}$, respectively.

 We will denote by $\cA_k$ the set of all indices $h$ as above, such that $D_h\in T(\ol{U}^{k-1})$ are exceptional divisors for $ f^{k-1}_k$. Note that all such $h$ differ only by permutations $g\in G$, for $G$ defined as in section $1$, and $\sum_{h\in \cA_k}D_h$ is the unique exceptional divisor in  $R(\ol{U}^{k-1})$. Let $\cA_{\geq k}$ be the set of all indices $h$, such that $D_h$ is an exceptional divisor for  $f^{k'-1}_{k'}$, and $k'\geq k$.



The structure of the forgetful morphism leads to a natural decomposition of classes in homology.

As a preliminary, we recall the following relation over 
$\ol{M}_{0,1}(\PP^n,d)$ (\cite{noi3}, Lemma  3.1 ) \bean 
 \psi^2-\sum_{h}N_{\psi h}  \psi D_h
+\sum_{h,h'}N_{hh'}  D_hD_{h'}+\frac{3}{d^4} k(H^2)^2-\frac{4}{d^3}k(H^3)=0  \eean where
$$ N_{\psi h} =\frac{|h|^2}{d^2}(6-4\frac{|h|}{d}), $$
$$ N_{h'h} =N_{hh'} =\left\{ \begin{array}{l} \frac{|h|^2}{d^2}(6\frac{|h'|}{d}-2\frac{|h|}{d}-3\frac{|h'|^2}{d^2}) \mbox{ if
} h'\supseteq h, \mbox{ and } \\
 -3\frac{|h|^2|h'|^2}{d^4}\mbox{ if
} h'\cap h=\emptyset. \end{array} \right.$$

With the notations from above, let $\psi(k)$ denote the relative cotangent cless for the morphism $f^{k}$, where $0\leq k\leq l_0$. In particular, $\psi(0)=\psi$, the $\psi$--class of $\ol{M}_{0,1}(\PP^n,d)$.

Relation (2.1) implies that $\psi(l_0)^2$ is the pullback of a class on $\ol{M}_{0,0}(X,d)$. In the case $d=2d'$, the same result is true for an appropriate definition of $\psi(l_0)$.

\begin{proposition}
Let $X$ be a projective variety. The following relation holds in $R(\ol{M}_{0,1}(X,d))$
\bea  \beta= -\frac{1}{2}[\psi f^*f_*\beta +f^*f_*(\psi\beta)-\sum_{k}j_{\bar{h} *}f_{\bar{h}}^*f_{h *}j_h^*\beta],\eea
where for each $k$ such that $0\leq k\leq l_0$, $h$ is an element in $\cA_k$.
\end{proposition}

\begin{proof}

The above formula is the decomposition of the class $ \beta $ due to successive blow-ups. We proceed by descending induction on $k$. The induction hypothesis states 
\bea \beta =  -\frac{1}{2}[\psi(k) f^{k *}f^{k}_*(\beta)+f^{k *} f^{k}_*(\psi(k)\beta)-\sum_{k'>k}( j_{\bar{h} *}f_{\bar{h}}^*f_{h *}j_h^*\beta+j_{h *}f_{h}^*f_{\bar{h} *}j_{\bar{h}}^*\beta )] \eea
where for each $k'$ such that $k< k'\leq l_0$, $h$ is an element in $\cA_{k'}$.

    When $k=l_0$, the relation on the projective bundle $U^{l_0}\to \ol{M}_{0,0}(X,d)$ is checked by applying  $\psi(l_0)^a$ to both sides, followed by $f^{l_0 }_*$ (the values $a=0, 1$ are enough). Indeed, we saw that $f^{l_0}_*(\psi(l_0)^2 \cap [U^{l_0}]^{vir}    ) =0$ while $f^{l_0}_*(\psi(l_0) \cap [U^{l_0}]^{vir} )=-2[ \ol{M}_{0,0}(X,d)]^{vir}$. When  $d=2d'$, the first induction step follows from the quisi-splitting of $U^{l_0}$ into a blow-up and a $\PP^1$--bundle.

 We now assume the formula true at the step $k$. By a blow-up calculation,
\bean \beta =  -j^{k-1 *}_{h} g^{*}_{h}  g_{h *}j^{k-1 *}_{h} \beta  + f^{k-1 *}_{k}  f^{k-1 }_{k *}\beta \eean
in $R(\ol{U}^{k-1})$.

  Apply the induction hypothesis to $ f^{k-1 }_{k *}\beta$ and substitute the result in equation (2.2).  Comparison of relative dualizing sheaves yields $f^{k-1 *}_{k}\psi(k) =  \psi(k-1)-\sum_{h\in\cA_k}D_h$. We recall that by Definition 1.4, $ \sum_{|h'|=|h|} D_{h'}\beta = j^{k-1}_{h *} j^{k-1 *}_{h}\beta$.

Furthermore,
$f^{k-1 *}j_{h\bar{h} *}=j^{k-1}_{\bar{h} *}f^{k-1 *}_{\bar{h}}+ j^{k-1}_{h *}f^{k-1 *}_{h}$
and therefore
 \bea & f^{k-1 *}f^{k-1}_*j^{k-1}_{h *}j^{k-1 *}_{h}\beta = f^{k-1 *}j_{h\bar{h} *}f^{k-1}_{h *} j^{k-1 *}_{h}\beta= & \\ &=j^{k-1}_{\bar{h} *}f^{k-1 *}_{\bar{h}}(f^{k-1}_{h *} j^{k-1 *}_{h}\beta) +  j^{k-1}_{h *}f^{k-1 *}_{h}(f^{k-1}_{h *} j^{k-1 *}_{h}\beta).& \eea
Also, $f^{k-1}_{k \bar{h}*}j^{k-1*}_{\bar{h}}+ f^{k-1}_{k {h}*}j^{k-1*}_{{h}}=j^{k-1 *}_{h\bar{h}}f^{k-1}_{k *}$, which implies
\bea  &j^{k-1}_{h *} j^{k-1 *}_{h}f^{k-1 *}f^{k-1}_*\beta =  j^{k-1}_{h *}f^{k-1 *}_{h}j^{*}_{h\bar{h}}f^{k-1}_*\beta =& \\ &=  j^{k-1}_{h *}f^{k-1 *}_{h}f^{k-1}_{\bar{h}*}j^{k-1*}_{\bar{h}}\beta  +  j^{k-1}_{h *}f^{k-1 *}_{h}  f^{k-1}_{{h}*}j^{k-1*}_{{h}}\beta.  &\eea
 Note that in fact $  f^{k-1 *}_{h}  f^{k-1}_{h *}=g^{*}_{h}  g_{h *} $.  Adding up  we obtain
\bea  &  f^{k-1 *}f^{k-1}_*j^{k-1}_{h *}j^{k-1 *}_{h}\beta+j^{k-1}_{h *} j^{k-1 *}_{h}f^{k-1 *}f^{k-1}_*\beta = &\\&= j^{k-1}_{\bar{h} *}f^{k-1 *}_{\bar{h}}(f^{k-1}_{h *} j^{k-1 *}_{h}\beta)+  j^{k-1}_{h *}f^{k-1 *}_{h}f^{k-1}_{\bar{h}*}j^{k-1*}_{\bar{h}}\beta + 2j^{k-1}_{h *}g^{*}_{h}  g_{h *}j^{k-1*}_{{h}}\beta.&\eea

Remark also that for any $h'\in \cA_{>k}$, $ j^{k *}_{h'}f^{k-1 }_{k *}\beta = f^{k-1 }_{k h' *}j^{k-1 *}_{h'}\beta$.

  Together all the above lead to the induction formula at step $k-1$.

\end{proof}

By the same logic, in the case $m\geq 1$ we obtain
\begin{proposition}
Let $X$ be a convex projective variety and $m\geq 2$. The following relation holds in $R(\ol{M}_{0,m}(X,d))$
\bea  \beta= \psi_1 f_m^*f_{m *}\beta +f_m^*f_{m *}(\beta\cdot D_{1,m})+\sum_{m\not\in h}f_{h\cup\{m\}}^*f_{h *}j_h^*\beta.\eea
\end{proposition}

Here $D_{1,m}$ denotes the Cartier divisor corresponding to the natural section $\sigma_1: \ol{M}_{0,1}(X,d)\to \ol{M}_{0,m}(X,d)$. The analogue to equation (2.1) in this context is
$$D_{1,m}^2+f_{m}^*\psi_1 \cdot D_{1,m}=0.$$

In particular we may apply the propositions above  to the tautological kappa classes, whose definition we recall below, following \cite{kk}.

\begin{definition} Let $f:\ol{M}_{0,n+1}(X,d)\to \ol{M}_{0,n}(X,d)$ denote the forgetful morphism and let $\alpha\in H^*(X)$. The kappa class  $$k_a(\alpha):=f_*(\psi^{a+1}ev^*\alpha  \cap [\ol{M}_{0,n+1}(X,d)]^{vir}).$$


 Assume $X\hookrightarrow \prod_{i=1}^s\PP^{n_i}$ such that any curve class $d\in H_2(X)$ may be assigned a unique tuple $(d_1,...,d_s)$, and let $\cD:=\bigsqcup_{i=1}^s\{1,...,d_i\}$. 
 For any $h=h_1\sqcup...\sqcup h_s\sqcup M_h \subset \cD\bigsqcup\{1,...,n\}$, let $|h|:=(|h_1|,...,|h_s|)$. Let $\bar{h}$ denote the complement of $h$ and consider the morphisms
\bea \diagram  
 \ol{M}_{0,M_h\cup\{ \star, \bullet\} }(X,|h|)\times_{X}\ol{M}_{0,M_{\bar{h}}\cup\{ \star\}}(X,|\bar{h}|)\dto^{f_h}
 \\ \ol{M}_{0,M_h\cup\{ \star\} }(X,|h|)\times_{X}\ol{M}_{0, M_{\bar{h}}\cup\{ \star\} }(X,|\bar{h}|)\rto  & \ol{M}_{0,n}(X,d), \enddiagram
  \eea
where $f_h$ is the map forgetting the marked point $\bullet$ on the first component, and identity on the second component. We define the class $$k_{h,a}(\alpha):= f_{h*}(\psi_{\bullet}^{a+1}ev_{\bullet}^*\alpha
\cap \Delta^{!}([\ol{M}_{0,l+2}(X,|h|)]^{vir}\otimes [\ol{M}_{0,n-l+1}(X,|\bar{h}|)]^{vir}) )$$
 in $H_*(\ol{M}_{0,l+1}(X,|h|)\times_{X}\ol{M}_{0,n-l+1}(X,|\bar{h}|)$ and in $T(\ol{M}_{0,n}(X,d))$. Here $\Delta:  \ol{M}_{0,l+2}(X,|h|)\times_{X}\ol{M}_{0,n-l+1}(X,|\bar{h}|)\to  \ol{M}_{0,l+2}(X,|h|)\times\ol{M}_{0,n-l+1}(X,|\bar{h}|)$ is the natural pullback of the diagonal $X\to X\times X$. 

When $a=-1$ we simply write $k(\alpha):=k_{-1}(\alpha)$, $k_h(\alpha):=k_{h,-1}(\alpha)$
\end{definition}

As a convention, for any class $\beta\in H^*(\ol{M}_{0,m}(X,d))$ we will denote  its image in $R(\ol{M}_{0,m}(X,d))$ by $\beta$ as well. 


\begin{corollary}
Let $X$ be a smooth projective variety and $\alpha\in H^*(X)$. 

(i) The following relation holds in $R(\ol{M}_{0,1}(X,d))$
\bea  ev_1^*\alpha = -\psi_1 k(\alpha) -k_0(\alpha)+k_{\Delta}(\alpha).\eea
Here $k_{\Delta}(\alpha):=\sum_{h}k_{h}(\alpha)$  in $T(\ol{M}_{0,1}(X,d))$, and each $k_{h}(\alpha)$ is the class on $\ol{M}_{0,1}(X,|h|)\times_{X}\ol{M}_{0,2}(X,|\bar{h}|)$ introduced in Definition 2.3. In other words, $k_{\Delta}(\alpha)$  is represented by the space of stable maps from split curves, of  one component intersecting the class $\alpha$ and another component containing the marked point.

(ii) Let $n\geq 2$. The following "change of variable" holds in $R(\ol{M}_{0,n}(X,d))$
\bea  ev^*_i\alpha-ev^*_j\alpha  =\psi_j k(\alpha) -k_{\Delta (i|j)}(\alpha | ).\eea
Here $k_{\Delta (i|j)}(\alpha | ):=\sum_{h}k_{h}(\alpha)$, where  $k_{h}(\alpha)$ is the  class on $\ol{M}_{0,\{ i, \star \} }(X,|h|)\times_{X}\ol{M}_{0, \{ j, \star \}}(X,|\bar{h}|)$ defined above.

\end{corollary}

\begin{proof}
Equation (i) follows from  Proposition 2.1, after observing that
 $$ f^*k_0(\alpha)=k_0(\alpha) -ev^*\alpha,$$
which is a special case of Lemma 5.6 in \cite{kk}.
For (ii), it is enough to pullback relation (i) via the two forgetful morphisms $\ol{M}_{0,2}(X,d)\to \ol{M}_{0,1}(X,d)$ and apply the change of variable formula \cite{leepan}, Theorem 1 for the $\psi$ classes.
Alternatively, we apply Proposition 2.2 and the comparison formula $f^*\psi_1=\psi_1-D_{1,m}$ for $\psi$ classes.
\end{proof}

 Let $H$ be a very ample divisor on $X$, let $|d|:=\int_d H$. Note that by Lemma 2.2 of \cite{pandharipande},
\bean k_0(\alpha)= -\frac{2}{|d|} k(H\alpha)+\frac{1}{|d|^2}k(\alpha)k(H^2)+\sum_{h}\frac{|h|^2}{|d|^2} k_{\bar{h}}(\alpha). \eean

 Equation (ii) in codimension one, together with the change of variable in $\psi$-- classes, has been employed by Lee and Pandharipande in \cite{leepan} for a reconstruction theorem in quantum cohomology, with an analogue in quantum $K$-theory.
 For codimension two classes in the two--plane Grassmaniann, this equation has also been proved previously in \cite{oprea2}. Recently Oprea informed us that he can also prove statement (ii) by other methods, involving projection formula via the forgetful map. Not surprisingly, the ensuing partial reconstruction for Gromov-Witten invariants (\cite{oprea2} for $G(2,n)$) may be obtained directly via WDVV equations. 

 Let $\cD=\bigsqcup_i \{1,...,d_i\}$. We will employ the notations introduced in Definition 1.3. $D_h$ will denote the class of $\ol{M}_h$ in  $T(\ol{M}_{0,m+1}(X,d))$.

\begin{theorem}
 Let $X$ be a smooth projective variety, $d\in H_2(X)$ a curve class and let $\alpha_1, ..., \alpha_s$ be generators for the ring $H^*(X)$. As an algebra over $T(\ol{M}_{0,m}(X,d))$, the extended tautological ring $T(\ol{M}_{0,m+1}(X,d))$ is generated by
\begin{itemize}
\item divisor classes $D_{1, m+1}$ and $\{D_{h}\}_{m+1\in h}$, where $h\subset \cD \bigsqcup \{2,...,m+1\}$ satisfies $h\not=\emptyset$ or $\{i\}$, $\forall i\in \{2,...,m+1\}$.
\item the image in $T(\ol{M}_{0,m+1}(X,d))$ of products $Q_h:=\prod_{i}k_{h}(\alpha_i)^{a_i}\in T(\ol{M}_{h})$ with $h$ as above. Here $a_i\geq 0$. 

\end{itemize}

\end{theorem}

Note that the first set of divisors may be further restricted due to the divisorial relations pulled back from $\ol{M}_{0,4}$ (see \cite{keel}). 

This theorem is complementary to \cite{noi3}, Theorem 3.4, where universal relations among the ring generators are given. In addition, specific relations in  $T(\ol{M}_{0,m+1}(X,d))$ come directly from relations in $H^*(X)$, induced by the kappa class decomposition
\bean k(\alpha\beta)=k(\alpha)ev_1^*\beta+k(\beta)ev_1^*\alpha+\psi_1k(\alpha)k(\beta)-\sum_{1\not\in h}k_h(\alpha)k_h(\beta). \eean 
 and its analogues on the boundary.

\begin{proof}
We will first prove a slightly different statement:
\begin{claim} Let $m> 0$. The algebra  $T(\ol{M}_{0,m+1}(X,d))$ is generated over $\QQ$ by $D_{1,i}$ and $ev^*_i\alpha$ for $i=1,...,m+1$, $k(\alpha)$,   $D_{h}$ and   $k_h(\alpha)$   (where $h$ does not necessarily contain $m+1$).

When $m=0$, $D_{1,i}$ is replaced by $\psi_1$ in the above.
\end{claim}
Here  $\alpha$ varies among all classes in $H^*(X)$. Note that by Lemma 2.2. of \cite{pandharipande} and Theorem 1 of \cite{leepan}, all the $\psi$-- classes at the marked points as well as nodes derive from the generators listed in the Claim. So are the classes $ev^*_h\alpha$ on the boundary, where $h$ is a node, due to Corollary 2.4 (ii).

 Let $f$ denote any of the forgetful maps between two of the moduli spaces or normal strata, let $g$ be any local embedding of boundary into a moduli space or into a boundary stratum, and $"\cdot"$  the tautological product on any stratum or moduli space. A monomial $\gamma$ in the tautological ring is obtained by any sequence of operations $f_*$, $g_*$, $"\cdot"$, out of an initial set of inputs $\{ev_j^*\alpha\}_j$, where $j$ is either a marked point or a node. The claim above follows independently of $m$, by induction on the number of operations $f_*$ performed inside the monomial. Indeed, when this number is $0$, the monomial is of the form $$\gamma=\prod_{l} ev^*_l\alpha_l \prod_{h}D_h^{a_h} \prod_{i,j}D_{i,j}^{b_{ij}},$$ in agreement with the claim. 
Let $M_{m+1}$ be the $\QQ$-vector space
$$< \prod_{l} ev^*_l\alpha_l \prod_{h}D_h^{a_h} \prod_{i,j}D_{i,j}^{b_{ij}}\prod_{r}k(\alpha_r)\prod_{r',h'}k_{h'}(\alpha_{r'}) >$$
in $T(\ol{M}_{0,m+1}(X,d))$. Note that $M_{m+1}$ is closed under $g_*$ and $\cdot$.
The induction step consists in checking that
$f_*(M_{m+1})\subseteq M_m$. Choose $\gamma\in M_{m+1}$ supported on a boundary stratum $\ol{M}_I$.
 Notice that on $\ol{M}_I$, $k(\alpha)=f^*k(\alpha)$ and  $k_h(\alpha)=f_h^*k_h(\alpha)$ if $m+1\not\in h$, while $k_{h}(\alpha)= k(\alpha)f^*D_{h\setminus\{m+1\}}- k_{h\setminus\{m+1\}}(\alpha)$ otherwise. After projection formula it remains to understand  the case $$ \gamma=\prod_{l} ev^*_l\alpha_l \prod_{h}D_h^{a_h} \prod_{i,j}D_{i,j}^{b_{ij}}.$$
Also by projection formula, the indices $l$ may be restricted to $l=m+1$, or nodes on components in the fiber which are contracted by $f$. We may also assume $j=m+1$. With the notations following Proposition 2.2,
 $$f_*(\beta D_{i,m+1}^b)=\sigma_i^*(\beta)\psi_i^{b-1}$$
for any class $\beta\in T^*(\ol{M}_{0,m+1}(X,d))$ and any integer $b\geq 1$. Assume $\gamma\not= 0$. Then for $l$ as above, $\sigma_i^*ev^*_l\alpha_l= ev_i^*\alpha_l$, while $\sigma_i^*D_h= D_h$ or $D_{h\setminus\{m+1\}}$, depending on whether $i\in h$ or not. 
 As all of the above are in $M_m$, it remains to study  $$\gamma= \prod_{l} ev^*_l\alpha_l \prod_{h}D_h^{a_h}.$$
We prove $f_*(\gamma)\in M_m$ by induction on the degrees $a_h$. If all $a_h\leq 1$, the statement follows from the additivity of $k$ classes (Lemma 3.3 in \cite{kk}). Assume $a_h\geq 2$ for some $h$ and assume that $m+1\not\in h$ (the other case works similarly).
 As $f^*D_{h}=D_{h\cup\{m+1\}}+D_h$, it follows that
 \bean D_h^{a_h}=D_h^{a_h-1}f^*D_{h}-D_h^{a_h-2}(D_{h\cup\{m+1\}}D_h).\eean
Furthermore, the support of $D_{h\cup\{m+1\}}D_h$  is isomorphic to the stratum $\ol{M}_{h}$ in $\ol{M}_{0,m}(X,d)$. The induced isomorphism on tautological rings identifies boundary with boundary, $ev_{m+1}^*\alpha$ with the pull-back at the node $ev_{h}^*\alpha$, and $D_h$ with $\psi_h$, the $\psi$-- class at the node of the generic curve represented by $\ol{M}_{h}$. These are all classes in $M_m$. Thus formula (2.5) provides the induction step that concludes the proof of the claim. 

 The claim implies the proposition after separating all the generators in  $f^*T(\ol{M}_{m}(X,d))$ from those of $T(\ol{M}_{m+1}(X,d))$. Note that $k(\alpha)=f^*k(\alpha)$  and $ev_i^*\alpha= f^*ev_i^*\alpha$ in $T(\ol{M}_{m+1}(X,d))$, for $i\not= m+1$. In the case $i=m+1$ we apply Corollary 2.4, (ii). 


 Finally, we reduce the $k$-classes $k(\alpha)$ to those of the ring generators $\alpha_0,..., \alpha_s$ via the decomposition formula (2.4).
This is again deduced from Corollary 2.4, (ii) for the marked points $m+1$ and $1$, after multiplication by $ev_{m+1}^*\beta$ and push-forward by the forgetful map.
We employ the comparison formula $f^*\psi_1=\psi_1-D_{1,m+1}$. 
 A similar relation on $\ol{M}_{0,0}(X,d)$ may be found in the next theorem.
\end{proof}

Recall the notations following Definition 1.3. Thus for $h\subset \cD$ and its complement $\bar{h}$, the space $\ol{M}_{h\bar{h}}$ is defined as $\ol{M}_{0,1}(X, |h|)\times_{X}\ol{M}_{0,1}(X, |\bar{h}|)$. Let $\pi_h$, $\pi_{\bar{h}}$ denote the projections on the components and let $\psi_h$, $\psi_{\bar{h}}$ be the pullbacks of the $\psi$-- class via the two projections. The sum $\psi_h+\psi_{\bar{h}}$ may be written as $-D_{h\bar{h}}^2$ in  $T(\ol{M}_{0,0}(X,d))$, where $D_{h\bar{h}}$ is the virtual fundamental class of  $\ol{M}_{h\bar{h}}$. We will denote by $ev_h=ev_{\bar{h}}$ the evaluation map at the node.

\begin{theorem}
 (i) Let $X$ be a smooth projective variety and $d\in H_2(X)$ a curve class. Let $\alpha_0, ..., \alpha_s$ be generators for the ring $H^*(X)$, such that $\alpha_0$ is the class of a very ample Cartier divisor.

Then the ring $T(\ol{M}_{0,0}(X,d))$ is generated by the following classes
\begin{itemize}
\item  $k(\alpha_0^3)$, $k(\alpha_i)$ and  $k(\alpha_0\alpha_i)$
for $i=1,...,s$;
\item the images in $T(\ol{M}_{0,0}(X,d))$ of  products $Q_h:=\prod_{i}k_{h}(\alpha_i)^{a_i}\psi_h^a \in T(\ol{M}_{h\bar{h}})$, where $a_i\geq 0$ are integers and $a=0$ or $1$.
\end{itemize}

(ii) Each relation $P(\alpha_0, ..., \alpha_s)=0$ in $H^*(X)$ induces relations $P_{-1}$, $P_{0}$, and $Q_h\{P_{h}\}_h$ in the extended tautological ring, where $Q_h$ are classes as above. Here $P_a:=k_a(P(\alpha_0, ..., \alpha_s))$ are polynomials of the generators listed above, derived by successive applications of the formulas
\begin{enumerate}
\item $ k(\alpha\beta)=  -\frac{1}{2}[k(\alpha) k_0(\beta)+ k_0(\alpha)k(\beta) - k_{\Delta}(\alpha |\beta ) ]. $

Here $ k_{\Delta}(\alpha |\beta ):=\sum_{h}k_{h}(\alpha)k_{\bar{h}}(\beta)$ on $\ol{M}_{0,1}(X,|h|)\times_{X}\ol{M}_{0,1}(X,|\bar{h}|)$. $k_0(\alpha)$ is given by  equation (2.3).
On the boundary

\item $k_{h}(\alpha\beta)=  -\frac{1}{2}[k_{h}(\alpha) k_{0}(\beta)+k_{h}(\beta) k_{0}(\alpha)+ (k_{h}(\alpha) k_{\bar{h}}(\beta)+k_{h}(\beta) k_{\bar{h}}(\alpha))\psi_{\bar{h}} -$\\$-k_{h}(\alpha)k_{\Delta_{\bar{h}}}(\beta) - k_{h}(\beta)k_{\Delta_{\bar{h}}}(\alpha)-\sum_{h'\subset h} (k_{h'}(\alpha)k_{h\setminus h'}(\beta )- k_{h'}(\beta)k_{h\setminus h'}(\alpha ) ) ].$

\end{enumerate}

\end{theorem}



\begin{proof}


By the previous claim,
elements in the extended tautological ring of $\ol{M}_{0,0}(X,d)$
 are come as products on various strata of classes $ ev_{h}^*\alpha$ and $f_*(ev^*\alpha\prod_{h}D_h^{a_h})$, with $a, a_h\geq 0$ (note that these include $k_h(\alpha)$).

If $a_h\leq 1$, then the last term is $k(\alpha)$ or a lienar combination of classes $k_{h'}(\alpha)$ on the boundary. If $a_h\geq 2$ for some $h$, then by the same procedure as in the proof of the claim we obtain products of classes $D_{h\bar{h}}$, $\psi_h$, $ ev_{h}^*\alpha$, or $k_h(\alpha)$ on the boundary.

 Moreover, Corollary 2.4, (i) applied to $\ol{M}_{0,1}(X,|h|)$ and $\ol{M}_{0,1}(X,|\bar{h}|)$, together with equations
\bea ev_h^*\alpha=ev_{\bar{h}}^*\alpha \mbox{ and }    k_{h,0}(\alpha)+k_{\bar{h},0}(\alpha)=k_0(\alpha) \mbox{ on } \ol{M}_{h\bar{h}} \eea
 result in the following formulas for  $ev_h^*\alpha$ and $ k_{h,0}(\alpha)$
\bea & ev_h^*\alpha  = -\frac{1}{2}[ k_0(\alpha)D_{h\bar{h}}+\psi_hk_{h}(\alpha)+ \psi_{\bar{h}} k_{\bar{h}}(\alpha)-k_{\Delta_{h}}(\alpha)-k_{\Delta_{\bar{h}}}(\alpha)],&\\
& k_{h,0}(\alpha)= \frac{1}{2}[ k_0(\alpha)D_{h\bar{h}}-\psi_hk_{h}(\alpha)+\psi_{\bar{h}} k_{\bar{h}}(\alpha)+k_{\Delta_{h}}(\alpha)-k_{\Delta_{\bar{h}}}(\alpha)].&\eea
Here $k_{\Delta_{h}}(\alpha):=\sum_{h'\subset h}k_{h'}(\alpha)D_{h\bar{h}}$. We note that $\psi_h+\psi_{\bar{h}}=-D_{h\bar{h}}$, $k_{h}(\alpha)+k_{\bar{h}}(\alpha)=k(\alpha)$ of $\ol{M}_{h\bar{h}}$. We have thus shown that the tautological ring of $\ol{M}_{h\bar{h}}$ is generated by classes $D_{h\bar{h}}$, $\psi_h$, and $k_{h}(\alpha)$, plus generators supported on the boundary of $\ol{M}_{h\bar{h}}$. Here $-D_{h\bar{h}}$ is identified to the first Chern class of the normal bundle $\cN_{\ol{M}_{h\bar{h}}|\ol{M}_{0,0}(X,d)}$. While multiplication with $D_{h\bar{h}}$ in $T(\ol{M}_{h\bar{h}})$ is in fact a product in $T(\ol{M}_{0,0}(X,d))$, classes of the form $\psi_h^a\prod_{\alpha}k_{h}(\alpha)^{a_{\alpha}}\in T(\ol{M}_{h\bar{h}})$ do not decompose in $T(\ol{M}_{0,0}(X,d))$. The only extra simplification we may afford is to assume $a\leq 1$, due to pullback of equation (3.1) from a product of projective spaces to $X$.


Consider generators $\alpha_0, ..., \alpha_s$ for the ring $H^*(X)$, and let $\alpha\in H^*(X)$. The second step is to decompose  $k(\alpha)$,  and $k_{h}(\alpha)$ into polynomials of $\{ k(\alpha_i)$, $k_{0}(\alpha_i)$ and $k_{h}(\alpha_i)\}_i$ with coefficients in $\QQ[\{ D_{h\bar{h}}, \psi_h\}_h]$. This is done by relations (1)-(2).

 Relation (1) is derived by multiplying both sides of equation (i), Corollary 2.4 by $ev^*\beta$, and pushing forward by the forgetful map $\ol{M}_{0,1}(X,d)\to \ol{M}_{0,0}(X,d)$. Note that $f^*k_0(\alpha)=k_0(\alpha)- ev_1^*(\alpha)$ as a special case of \cite{kk}, Lemma 5.6.


To justify relation (2), we combine the analog of equation (2.4) on $\ol{M}_{h\bar{h}}$
 \bea k_{h}(\alpha\beta)=  k_{h}(\alpha) ev_{h}^*(\beta)+ ev_{h}^*(\alpha)k_{h}(\beta) 
+\psi_hk_h(\alpha)k_h(\beta)+  \sum_{h'\subset h} k_{h'}(\alpha) k_{h'}(\beta ) ], \eea
  with the formula for $ev_h^*(\beta)$, written above.

 We note that in the polynomial ring over $\QQ$ generated by the  ring generators listed above,  polynomials $P_a$ with $a>0$ are in the ideal generated by $P_0$, $P_{-1}$, and the boundary $Q_hP_{h, -1}$, due to equation (2.1).



\end{proof}


\begin{remark}
The generators listed in the two previous theorems take an especially simple form when the cohomology of $X$ is generated by divisor classes. $T(\ol{M}_{0,0}(X,d))$ is generated by $k(\alpha_0^3)$, $k(\alpha_0\alpha_i)$, and boundary classes $D_{h\bar{h}}$ and $F_h$, the image of $\psi_h$ in  $T(\ol{M}_{0,0}(X,d))$. $T(\ol{M}_{0,m+1}(X,d))$ is generated over $T(\ol{M}_{0,m}(X,d))$ by boundary classes  $D_h$  and $D_{1,m+1}$ ($\psi_1$ in the case $m=0$). \end{remark}


\section{Maps to $\PP^n$}

In the case of maps to $\PP^n$,it turns out that the tautological relations described in Theorem 2.6, together with a set of natural universal relations on the boundary, completely determine the tautological ring structure. We recall that in this case  the moduli space is a smooth Deligne-Mumford stack and its entire  cohomology is tautological. The universal relations mentioned above are described in the following lemma.   

\begin{notation}

  Let $\ol{M}_{h\bar{h}}:=\ol{M}_{0,1}(\PP^n,|h|)\times_{\PP^n}\ol{M}_{0,1}(\PP^n,|\bar{h}|)$. Consider the boundary map
$g_{h\bar{h}}: \ol{M}_{h\bar{h}} \to
\ol{M}_{0,0}(\PP^n,d)$ and its class $D_{h\bar{h}}$ in
$T(\ol{M}_{0,0}(\PP^n,d))$. Let $\pi_h$ and $\pi_{\bar{h}}$ denote
the projections on the two components and let
$\psi_h:=\pi_h^*\psi_1$,    $\psi_{\bar{h}}:=\pi_{\bar{h}}^*(\psi_1)$. The images of the classes  $\psi_h$ and $\psi_{\bar{h}}$ in
$T(\ol{M}_{0,0}(\PP^n,d))$ will be denoted by $F_h$ and $F_{\bar{h}}$.
\end{notation}

\begin{notation}

Let $H$ be the hyperplane divisor in $\PP^n$. Define 
\bea &  b_h:= -N_{\psi\bar{h}}D_{h\bar{h}}+\sum_{h'\subset
h}N_{\psi h'} D_{h'\bar{h}'} - \sum_{h'\subset
\bar{h}}N_{\psi h'} D_{h'\bar{h}'}, & \\ &
c_h:=\frac{3}{d^4}k(H^2)^2-\frac{4}{d^3}k(H^3) + \sum_{(h'\subset h) \lor
(h'\subset\bar{h})}N_{h'h'}  D_{h'\bar{h}'}^2 +&\\ &  +\sum_{(h''\not= h'\subset h) \lor
(h''\not= h'\subset\bar{h})}N_{h''h'}  D_{h'\bar{h}'}
D_{h''\bar{h}''}  +N_{\bar{h}\bar{h}} D_{h\bar{h}}^2- &\\ & -2\sum_{h'\subset
h}N_{h'\bar{h}} D_{h\bar{h}}D_{h'\bar{h}'} +2\sum_{h'\subset
\bar{h}}N_{h'\bar{h}} D_{h\bar{h}}D_{h'\bar{h}'}
& \eea where $N_{\psi h}$ and $ N_{hh'}$ are as in equation (2.1).
\end{notation}

\begin{lemma}
The algebra $T(\ol{M}_{0,0}(\PP^n,d))$ is generated over
$H^*(\ol{M}_{0,0}(\PP^n,d))$ by codimension one classes
$\{D_{h\bar{h}}\}_{h\subset D}$ and codimension two classes
$\{F_h\}_{h\subset D}$, satisfying the following relations:
\begin{enumerate}
\item $D_{h\bar{h}}D_{h'\bar{h}'}=0$ whenever the set $\{h\bigcap h', h\bigcap \bar{h}', \bar{h}\bigcap h', \bar{h}\bigcap\bar{h}'\}$ has four distinct elements;
\item $F_h+F_{\bar{h}}=-D_{h\bar{h}}^2$;
\item $F_hD_{h'\bar{h}'}+F_{h'}D_{h\bar{h}}=\sum_{h''; \bar{h}'\subset h''\subset h}D_{h\bar{h}}D_{h'\bar{h}'}D_{h''\bar{h}''}$ for any $h, h'\subset D$ such that $h\cap h'\not=\emptyset$;
\item $ F_h^2 +b_h F_hD_{h\bar{h}} +c_h D_{h\bar{h}}^2= 0$
\item $F_hF_{h'} = (b_hD_{h\bar{h}}+\sum_{h''; \bar{h}'\subset h''\subset h}D_{h''\bar{h}''} )D_{h'\bar{h}'}F_h+c_hD_{h\bar{h}}^2$ for any $h, h'\subset D$ such that $h\cap h'\not=\emptyset$.
\end{enumerate}

\end{lemma}

\begin{proof}

By definition, $T(\ol{M}_{0,0}(\PP^n,d))$ is the extended cohomology
ring of the $S_d$--network generated by the regular local embeddings
$ g_{h\bar{h}}$ ( see \cite{noi1}).  Thus the
generators of $H^*(\ol{M}_{h\bar{h}})$ for all $h\in \{1,...,d\}$
form a complete set of generators for $T(\ol{M}_{0,0}(\PP^n,d))$
over $H^*(\ol{M}_{0,0}(\PP^n,d))$.

Let $f_h: \ol{M}_{0,1}(\PP^n,|h|) \to \ol{M}_{0,0}(\PP^n,|h|)$,
$f_{\bar{h}}: \ol{M}_{0,1}(\PP^n,|\bar{h}|) \to
\ol{M}_{0,0}(\PP^n,|\bar{h}|) $ be the forgetful morphisms. Denote by $H_h:=ev_h^*H$, where $ev_h:  \ol{M}_{0,1}(\PP^n,|h|) \to \PP^n$ is the usual evaluation map, and similarly for $\bar{h}$. 

  By Theorem 3.23 in \cite{noi1}, the algebra generators of
$H^*(\ol{M}_{0,1}(\PP^n,|h|))$ over $\QQ$ are 
$\psi_h$, $H_h$ and boundary divisors $D_{h'\bar{h}'}$ for
$h'\subset h$. Similarly, the generators of
$H^*(\ol{M}_{0,1}(\PP^n,|\bar{h}|))$ are 
$\psi_{\bar{h}}$, $H_{\bar{h}}$ and boundary divisors $D_{h'\bar{h}'}$ for
$h'\subset {\bar{h}}$. By Lemma 2.2, \cite{pandharipande}, $k_h(H^2)$
can replace $H_h$ among the generators. Moreover, these two divisors
can also be written as linear combinations of $k(H^2)$, $\psi_h$ and
boundary divisors, via Lemma 2.2, \cite{pandharipande} applied to both $\ol{M}_{0,1}(\PP^n,|h|)$ and $\ol{M}_{0,1}(\PP^n,|\bar{h}|)$, and the relations $H_h=H_{\bar{h}}$, $k_h(H^2)+k_{\bar{h}}(H^2)=k(H^2)$ and $\psi_h+\psi_{\bar{h}}=-D_{h\bar{h}}$ on $\ol{M}_{h\bar{h}}$. Thus\bea&  k_h(H^2)=|h||\bar{h}|\psi_h-\frac{|h||\bar{h}|^2}{d}D_{h\bar{h}}+ & \\ &+ \frac{|h|}{d}k(H^2)+\sum_{h'\subset h} \frac{|h'|^2|\bar{h}|}{d} D_{h'\bar{h}'}-  \sum_{h'\subset \bar{h}} \frac{|h'|^2|h|}{d} D_{h'\bar{h}'}. & \eea
 Following the above analysis, each $H^*(\ol{M}_{h\bar{h}})$ contributes exactly two new generators $D_{h\bar{h}}:=[\ol{M}_{h\bar{h}}]$ and $F_h:=\psi_hD_{h\bar{h}}$ to the algebra $T(\ol{M}_{0,0}(\PP^n,d))$ over $H^*(\ol{M}_{0,0}(\PP^n,d))$. A priory, one should also consider classes  $\psi_h^kD_{h\bar{h}}$ with $k>1$, but these will be written in terms of the classes above, due to the existence of a quadratic relation in $\psi_h$ on $\ol{M}_{h\bar{h}}$.
Indeed, summing equation (2.1) for $\ol{M}_{0,1}(\PP^n,|h|)$ and $\ol{M}_{0,1}(\PP^n,|\bar{h}|)$ yields $k(H^3)$ in terms of $\psi_h$, $\psi_{\bar{h}}$, $k_h(H^2)$, $k_{\bar{h}}(H^2)$ and boundary divisors. We have seen that $\psi_{\bar{h}}$, $k_h(H^2)$, $k_{\bar{h}}(H^2)$ are all expressions of $k(H^2)$, $\psi_h$ and boundary divisors. The following relation is ensues on $\ol{M}_{h\bar{h}}$:
\bean  \psi_h^2 +b_h \psi_h +c_h D_{h\bar{h}}= 0. \eean
 Equation (2) is a standard relation between $\psi$--classes. Equation (3) is a corollary of Theorem 1 in \cite{leepan}. 
  Relation (4) follows from equation (3.1) and, in conjunction with Theorem 1, \cite{leepan}, implies relation (5).

\end{proof}

We will keep the notations from section 2 throughout.
When the target is $\PP^n$, the factorization of the forgetful map introduced in section 2 is $$\ol{M}_{0,1}(\PP^n,d)=\ol{U}^0\to...\to \ol{U}^{\lfloor (d-1)/2\rfloor }\to \ol{M}_{0,0}(\PP^n,d),$$
such that $D_{\bar{h}}$ is the exceptional divisor of $\ol{U}^{|h|-1}$ whenever $|h|\leq (d-1)/2$, while the support of  $D_h$ is the strict transform of a section $S_{\bar{h}} \subset \ol{U}^{|h|}_{h\bar{h}}\hookrightarrow \ol{U}^{|h|}$. We will also denote by $S_{\bar{h}}$ the image of this section in  $\ol{U}^{k}$, when $k>|h|$. The structure of $T(\ol{M}_{0,1}(\PP^n,d))$ as an algebra over $T( \ol{M}_{0,0}(\PP^n,d))$ is completely determined due to the sequence above: each divisor $D_{\bar{h}}$ is annihilated by $\Ker( T(\ol{U}^{|h|-1})\to T(S_{\bar{h}})$, and satisfies a quadratic equation (\cite{noi3}, Theorem 3.3). 

The following corollary lists low degree relations in $T(\ol{M}_{0,1}(\PP^n,d))$. These have been proven already in \cite{noi1}, Theorem 3.23. However, our point here is that these relations derive from Lemma 3.1. due to the structure of the extension $T( \ol{M}_{0,0}(\PP^n,d))\to T(\ol{M}_{0,1}(\PP^n,d))$ described above. Indeed, Relation (1) is derived from Lemma 3.1 (1), while relation (2) is derived from Lemma 3.1 (3).

\begin{corollary}
The following relations hold in $T(\ol{M}_{0,1}(\PP^n, d))$.
\begin{enumerate}
\item $D_hD_{h'}=0$ whenever $h\cap h'\not= h, h'$ or $\emptyset$.
\item $D_hD_{h'}(\psi+\sum_{h''\subseteq h\cup h'} D_{h''})$ for all $h, h'$ such that $h\cap h'=\emptyset$. \end{enumerate}

\end{corollary}

\begin{proof}
 Assume $0<|h'|\leq |h|<\lfloor  d/2\rfloor$. If $D_{h\bar{h}}D_{h'\bar{h}'}=0$, then on any intermediate space $\ol{U}^k$ with $|h|\leq k\leq \lfloor  d/2\rfloor$, restriction to the section $S_{\bar{h}}$ yields $s_{\bar{h}}D_{h'\bar{h}'}=0$, where
$s_{\bar{h}}$ denotes the class of $S_{\bar{h}}$ on $\ol{U}^{|h|}_{h\bar{h}}$. Thus by the considerations above, $D_{\bar{h}}D_{h'\bar{h}'}=0$ and also $D_{h}D_{h'\bar{h}'}=0$ on all spaces $\ol{U}^k$ with $|h|\geq k$. After restriction to $S_{\bar{h}'}$ and reiteration of the argument, we obtain   $D_{h}D_{\bar{h}'}=0$, $D_{\bar{h}}D_{\bar{h}'}=0$, and thus also $D_{\bar{h}}D_{h'}=0$ and $D_{h}D_{h'}=0$ on all spaces $\ol{U}^k$ with $|h'|\geq k$. We recall that $D_{h}+D_{\bar{h}}=D_{h\bar{h}}$.

 We now prove relation (2) under the same assumption as above.
  The other cases work similarly. First we notice the following 
\begin{claim}
$ D_{h\bar{h}}( a+ b s_{\bar{h}})=0 \mbox{ on } \ol{U}^{|h|}
\Rightarrow  D_{h}( a- b D_{\bar{h}})=0\mbox{ on }
\ol{U}^{|h|-1}$ whenever $0< |h|<\lfloor  d/2\rfloor $.
\end{claim}

Indeed, the quadratic relation satisfied by the exceptional divisor $D_{\bar{h}}$ may be written as
$$(D_{h\bar{h}}-D_{\bar{h}})( s_{\bar{h}}-D_{\bar{h}})=0.$$
This immediately implies the claim. 

Recall that any  $\psi$-- class  on $\ol{M}_{0,3}(\PP^n,d)$ may be written in terms of boundary divisors. (This is, for example, a consequence of \cite{leepan}, Theorem 1). Let $\ol{U}^{|h|}_{h, \bar{h}\bigcap \bar{h'}, h'}$ be the normal stratum of $\ol{U}^{|h|}$ which is birational to $\ol{M}_{0,1}(\PP^n, |h|)\times_{\PP^n} \ol{M}_{0,3}(\PP^n, |\bar{h}\bigcap \bar{h'}|)\times_{\PP^n}\ol{M}_{0,1}(\PP^n, |h'|)$. The following analogue of the above property holds here
\bean D_{h\bar{h}}D_{h'\bar{h'}}[ \psi(|h|)+s_{\bar{h}}+s_{\bar{h'}}-\sum_{h''\supseteq
h\cup h', |h''|<|\bar{h}|}D_{h''} ]=0,\eean
 where $\psi(|h|)$ is the first Chern class of the
 relative dualizing sheaf for $f^{|h|}:\ol{U}^{|h|}
 \to\ol{M}_{0,0}(X,d)$. Therefore, its pullback to $\ol{U}^{|h|}_{h, \bar{h}\bigcap \bar{h'}, h'}$ differs from the usual $\psi$--class by $s_{\bar{h}}+s_{\bar{h'}}$.

By the Claim above, this equation implies
$$D_{h}D_{h'\bar{h'}}(\psi(|h'|)+s_{\bar{h'}}-\sum_{h''\supseteq
h\cup h', |h''|<|\bar{h'}|}D_{h''})=0$$ on $\ol{U}^{|h'|} $. Indeed,
$\psi(|h'|)=f^{|h'| *}_{|h|}\psi(|h|)+\sum_{|\bar{h}|\leq
|h''|<|\bar{h'}|}D_{h''}$ and $s_{\bar{h'}}=f^{|h'|
*}_{|h|}s_{\bar{h'}} -\sum_{|\bar{h}|< |h''|, h''\subset \bar{h'}}D_{h''}$ on
$U^{|h'|}_{h'\bar{h'}}$, as $s_{\bar{h''}}\subset s_{\bar{h'}}$ on
$U^{|h''|}_{h'\bar{h'}}$ for all $h''\subset \bar{h'}$. Relation (1) also intervened in the computation. Again by the Claim, the desired relation (2) is obtained. 

We note that equation (3.2) is the pullback of relation (3) in Lemma 3.1. Indeed, for any integer $k$ such that $|h|\leq k<\lfloor d/2\rfloor$,
$$ s_{\bar{h}}=(f^{k *}\psi_{\bar{h}}+\psi(k)+\sum_{h\subset h''\in I, |h''|<|\bar{h}|}D_{h''}) \mbox{ on } \ol{U}^k_{h\bar{h}},$$
and an analogue relation exists for $s_{\bar{h}'}$ (\cite{noi3}, Theorem 3.3). Adding up these two formulas leads to equation (3.2).

\end{proof}

\begin{notation}
Consider the vector $V_{l}$ with entries in $T(\ol{M}_{0,0}(\PP^n,d))$
$$V_{l}:=(k(H^{l+1}), k(H^{l}), \{k_{h}(H^{l}) \}_{h\in I}, \{\psi_hk_{h}(H^{l}) \}_{h\in I})^T.$$
\end{notation}
Theorem 2.6 gives a formula for the transition matrix $A$ with entries in $\QQ[k(H^2), k(H^3), \{D_{h\bar{h}}, F_h\}_{h\in I}]$, such that $V_{l+1}=AV_l$ for all $l\geq 1$. Indeed, specializing equation (1) of Theorem 2.6 for $\alpha=H^{l-1}, \beta=H^2$, and equation (2) for $\alpha= H^l$, $\beta =H$ yields all the elements of $A$. Equation (2.3) is necessary for the transition matrix calculation, together with the formula for $k_{h}(H^2)$ derived in the proof of Lemma 3.1.

 The open stratum component $A^o$ of $A$ will be of special interest. It is the transition matrix for the vector $V^o_l:=(k(H^{l+1}), k(H^{l}))^T$ in $H^*(M_{0,0}(\PP^n,d))$.
\bean A^o=\left( \begin{array}{cc} \frac{1}{d}k(H^2) & \frac{1}{d}k(H^3)-\frac{1}{d^2}k(H^2)^2 \\ 1 & 0 \end{array}\right),\eean
while $V^o_1 =\left(\begin{array}{c} k(H^2)\\d\end{array}\right)$ and $V^o_{l+1}=(A^{o})^lV^o_1$ for $l\geq 1$.

\begin{theorem}
The extended cohomology ring $T(\ol{M}_{0,0}(\PP^n,d))$ is the $\QQ$-- algebra generated by classes $k(H^2)$, $k(H^3)$, and  $\{D_{h\bar{h}}, F_h\}_{h\in I}$. A complete set of relations consists of all relations in Lemma 3.1, together with the vectorial relation $V_{n+1}=0$.
\end{theorem}

\begin{proof}
 The ring generators of $T(\ol{M}_{0,0}(\PP^n,d))$ have already been determined in Theorem 2.6. Indeed, we recall that $H^*(\PP^n;\QQ)=\QQ[H]/(H^{n+1})$.

Let $f:\ol{M}_{0,1}(\PP^n,d)\to \ol{M}_{0,0}(\PP^n,d)$ denote the forgetful map.
In \cite{noi3}, $T(\ol{M}_{0,1}(\PP^n,d))$ is presented as an algebra over $T(\ol{M}_{0,0}(\PP^n,d))$, via the intermediate extensions $\{T(\ol{U}^l)\}_{1\leq l< d/2}$. Another presentation for the ring $T(\ol{M}_{0,1}(\PP^n,d))$ appears in \cite{noi1}.  We recall here that the ring generators of $T(\ol{M}_{0,1}(\PP^n,d))$ are $H$, $\{D_h\}_{h\subset\{1,...,d\} }$, and $k(H^2)=f^*k(H^2)$ (or, alternatively, $\psi$). A complete set of relations consists of those listed in Corollary 3.2, plus $ev^*H^{n+1}=0$, $P=0$ and $D_hP_h=0$, where $P$, $P_h$ are degree $n$ -- polynomials of variables $ev^*H$, $\{D_h\}_{h\subset\{1,...,d\} }$, and $f^*k(H^2)$, both involving the monomial $f^*k(H^2)^n$.

 In fact, $P$ is the polynomial relation of minimal degree having a summand of the type $c\cdot f^*k(H^2)^l$, where $c$ and $l$ are positive constants. On the other hand, $f^*k(H^{n+1})$ satisfies the same property, as shown by formula (3.3). This forces the equality $P= af^*k(H^{n+1})$ modulo relations (1) and (2) of Corollary 3.2; here $a$ is a constant. Similarly, modulo $P$, $P_h$ is the minimal degree polynomial relation having a summand of the type $c\cdot D_hf^*k(H^2)^l$. By formula (3.3) for the open stratum of $\ol{M}_{h\bar{h}}$, and due to the dependence of $k_h(H^2)$ on $k(H^2)$ expressed in the Proof of Lemma 3.1, we again obtain that $P_h= a_hf^*k_{h}(H^{n+1})$ modulo the low degree relations on the boundary mentioned above. Here $a_h$ is a constant. Furthermore, by Corollary 2.4, (i)
\bea  ev^*H^{n+1} = -\psi k(H^{n+1}) -k_0(H^{n+1})+k_{\Delta}(H^{n+1}).      \eea
 By equation (2.3), $k_0(H^{n+1})$ depends linearly on $k(H^{n+2})$, $k_{h}(H^{n+1})$ and $k(H^{n+1})k(H^2)$.

 Let $R$ be the $\QQ$-- algebra generated by $k(H^2)$, $k(H^3)$, $\{D_{h\bar{h}}, F_h\}_{h\in I}$, with the ideal of relations described in the statement of the present theorem. By the considerations above, in conjunction with Lemma 3.1, we conclude that
   $$ T(\ol{M}_{0,1}(\PP^n,d))\cong R[\psi, \{D_h\}_h]/\cJ, $$
where  $D_h+D_{\bar{h}}=D_{h\bar{h}}$ and the ideal $\cJ$ consists of relation in $ T(\ol{M}_{0,1}(\PP^n,d))$ over $ T(\ol{M}_{0,0}(\PP^n,d))$ (Theorem 3.3, \cite{noi3}). Note again that $\psi$ may be substituted by $H$ via Lemma 2.2 in \cite{pandharipande}. The structure of $\cJ$ implies 
$$T(\ol{M}_{0,0}(\PP^n,d))\cong R.$$
 Indeed, the $R$-- algebra $T(\ol{M}_{0,1}(\PP^n,d))$ is constructed by a sequence of algebra extensions of the type
$R_{i+1}=R_{i}[x]/(x^2+bx+c, Vx)$, where $b,c\in R_i$, and $V$ is a vector with entries in $R_i$. The only relations in $R_{i+1}$ which are not relations in $R_i$, but involve only elements in $R_i$ are multiples of $Vx(x+b)=-Vc$. Corollary 3.2 shows that in the present case, all such relations  are already justified by equations (1)-(3) of Lemma 3.1.
\end{proof}

$H^*(\ol{M}_{0,0}(\PP^n,d))$ is the $S_d$-invariant subalgebra of $T(\ol{M}_{0,0}(\PP^n,d))$. Listing ring  generators and relations for $H^*(\ol{M}_{0,0}(\PP^n,d))$ is difficult due to the intricate combinatorics of stable trees involved in the structure of the boundary. On the other hand, an additive basis for $H^*(\ol{M}_{0,0}(\PP^n,d))$ in terms of decorated trees, may be extracted from Theorem 3.3, as done in \cite{noi2} for the case when the number of marked points was positive. Thus $T(\ol{M}_{0,0}(\PP^n,d))$ is a rather compact, geometrically suggestive way of encoding the cohomology of $\ol{M}_{0,0}(\PP^n,d)$. Aside from the problem of tree combinatorics, Lemma 3.1 and formula (3.1) can simplify the task of extracting the ring structure of $H^*(\ol{M}_{0,0}(\PP^n,d))$ from $T(\ol{M}_{0,0}(\PP^n,d))$. The following lemma illustrates this idea.
 Let $s:T(\ol{M}_{0,0}(\PP^n,d))\to H^*(\ol{M}_{0,0}(\PP^n,d))$ be the Reynolds operator:
for any class $\gamma$, $s(\gamma):=\sum_{\sigma\in S_d}\sigma_*(\gamma)$.

\begin{lemma}
 Let $\alpha_h$ denote a class in $H^*(\ol{M}_{h\bar{h}})$, as well as its image in $T(\ol{M}_{0,0}(\PP^n,d))$. The following ideals are equal in $H^*(\ol{M}_{0,0}(\PP^n,d))$
\bea & H^*(\ol{M}_{0,0}(\PP^n,d))\bigcap (\{\alpha_{\sigma (h)}\}_{\sigma\in S_d}) = & \\ &=( \{ s(\alpha_{h}D_{I}),  s(\psi_{h}\alpha_{h}D_{I}),  s(\alpha_{h}D_{h\bar{h}}D_{I})\}_{I} ),&\eea
where all choices of $I=\{h_1, ..., h_l\}$ are considered with the properties that $h\not\in I$, $\{h'\bar{h}'\}\not\subset I$ for any $h'$, and $D_I:=\prod_{h'\in I}D_{h'\bar{h}'}^{b_{h'}}$ is a nonzero boundary class in
$T(\ol{M}_{0,1}(\PP^n,d))$.
\end{lemma}
\begin{proof}
An element in the ideal $H^*(\ol{M}_{0,0}(\PP^n,d))\bigcap (\{\alpha_{\sigma (h)}\}_{\sigma\in S_d})$ is of the form $ \sum_{\sigma\in S_d}\alpha_{\sigma (h)}\beta_{\sigma (h)}\beta_{\sigma (I)}$, where $I$ is as above, $\beta_I:=\prod_{h'\in I}\beta_{h'}$ for classes $\beta_{h'}$ supported on $\ol{M}_{h'\bar{h}'}$.  Theorem 3.3 implies that it is enough to consider $\beta_{h'}=F_{h'}^aD_{h'\bar{h}'}^b$ for some non-negative integers $a$ and $b$. Assume $h\bigcap h'\not=\emptyset$. If $a>0$, then by Lemma 3.1, (3)  $$\alpha_{h} F_{h'}=-\alpha_{h}\psi_{h} D_{h'\bar{h}'}+\alpha_{h}\sum_{h'\subset h''\subset \bar{h}} D_{h''\bar{h}''}.$$ Thus we may assume $a=0$.

 Similarly, $\beta_{h}=\psi_{h}^aD_{h\bar{h}}^b$ for some non-negative integers $a$ and $b$. Due to formula (3.1), we may consider $a=0$ or $1$.  In the former case, assume $b\geq 1$ and let $\gamma_h:=D_{h\bar{h}}^{b-1}D_I$.  Then the following equality
\bea  s(D_{h\bar{h}}^bD_I\alpha_h)=s(\alpha_hD_{h\bar{h}})s(\gamma_h)-s(\alpha_hD_{h\bar{h}}\gamma_{h'}) \eea
permits the inductive decrease of $b$ down to $b=1$.

   If $a=1$, let $\gamma_h:=\alpha_hD_{h\bar{h}}^{b-1}D_I$.  Then via equation Lemma 3.1,(3)
\bea s(D_{h\bar{h}}\psi_h\gamma_h)=\frac{1}{2}[s(F_h)s(\gamma_h)+s(\gamma_h\psi_h)s(D_{h\bar{h}})-s(D_{h'\bar{h}'}D_{h''\bar{h}''}\gamma_h)], \eea
where $h'\not=h$ is such that $|h'|=|h|$, $h'\bigcap h=\emptyset$, and $h'\subset h''\subset \bar{h}$. Inductively, this formula permits the decrease of $b$ down to $b=0$.
\end{proof}

\begin{example}
The above lemma is rather imprecise in the set of choices for $D_I$. For example, when $d=3$, the boundary combinatorics is simple, and the lemma may be recast as
\bea & H^*(\ol{M}_{0,0}(\PP^n,3))\bigcap (\{\alpha_{\sigma (h)}\}_{\sigma\in S_3}) = & \\ &=( s(\alpha_{h}),  s(\psi_{h}\alpha_{h}),  s(\alpha_{h}D_{h\bar{h}}), s(\alpha_{h}D_{h\bar{h}}^2)  ),&\eea
or, alternatively, the last two terms may be replaced by  $s(\alpha_{h}D_{h'\bar{h}'}D_{h''\bar{h}''})$, $s(\alpha_{h}D_{h'\bar{h}'}).$
Indeed, keeping the same notations as in the proof of the lemma, and with the convention $|h|=|h'|=|h''|=2$, one simplifies
$$s(\alpha_hD_{h'\bar{h}'})=s(\alpha_h)s(D_{h\bar{h}})-s(\alpha_hD_{h\bar{h}}),$$
etc. The ring structure of $H^*(\ol{M}_{0,0}(\PP^n,3))$ follows.

Let $P_2$ denote the set whose elements are the three cardinal 2 subsets of $\{1,2,3\}$.
 The unique boundary divisor in  $H^*(\ol{M}_{0,0}(\PP^n,3))$  is  $\sigma_1:= \sum_{h\in P_2}D_{h\bar{h}}$. The unique codimension two substratum in $H^*(\ol{M}_{0,0}(\PP^n,3))$ is $s_2:= \sum_{h\not=h'\in P_2}D_{h\bar{h}}D_{h'\bar{h}'}$ and the unique codimension three substratum is $\sigma_3:= \prod_{h\in P_2}D_{h\bar{h}}$. For consistency with \cite{behrend2}, let $\sigma_2:= \sigma^2_1-4s_2$. Define $\tau:=\frac{1}{2}\sum_{h=1}^3F_{\bar{h}}$ and
$$ \rho:=\frac{1}{4}(\psi')^2=-\frac{1}{27}[7\tau +\frac{1}{4}\sigma_1^2+\frac{1}{2}\sigma_2 +\frac{1}{4} k(H^2)^2 - k(H^3)] $$
 conform equation (2.1). Here $\psi'=\psi-\sigma_1$.
\begin{corollary} The algebra  $H^*(\ol{M}_{0,0}(\PP^n,3))$ is generated  over $\QQ$  by classes $k(H^2)$, $\rho$, $\sigma_1$, $\sigma_2$, $\sigma_3$ and $\tau$.
 The ideal of relations is $$\cI=( \tau\sigma_3, \rho\sigma_3, \tau^2-\rho\sigma_2, U_{n+1} ).$$
\end{corollary}

Here $ U_{l}$ is the vector \bea U_l=\left(\begin{array}{c} \sum_{h=1}^3 D_{h\bar{h}}k_{h}(H^{l})\\ \sum_{h=1}^3 \psi_hk_{h}(H^{l}) \\     \sum_{h=1}^3k_{h}(H^{l})  \\  k_0(H^{l}) \\ k(H^{l}) \end{array}\right)\eea

Modulo an elementary transformation of the vector $U_l$, this corollary is Conjecture 4.19 in \cite{behrend2}. Our particular choice of classes $\rho$, $\sigma_2$ was designed to fit the form of the conjecture.

 The product $\psi_h k_{h}(H^{l})$ is the usual cup product in $H^*(\ol{M}_{h\bar{h}})$, and the result is thought of as a class in $T(\ol{M}_{0,0}(\PP^n,3))$. The recursion matrix for $\{U_l\}_l$ is derived from the transition matrix $A$ in a straightforward way.
\begin{proof}
The list of generators as well as relation $U_{n+1}$ are consequences of the previous Lemma (in its more precise version above). Notice that this would include $s(k_{h}(H^{n+1})D_{h'\bar{h}'}D_{h''\bar{h}''})$ among relations. By additivity of the kappa classes though, this coincides with $k(H^{n+1})\sigma_3$, relation implied by $U_{n+1}$.

The fixed relations $ \tau\sigma_3, \rho\sigma_3, \tau^2-\rho\sigma_2$ are derived from Lemma 3.1 and equation (3.1). Indeed, for any degree $d$,  special elements  in the ideal of relations (1) and (3) of Lemma 3.1 are
$$  \prod_{i=1}^sD_{h_i\bar{h}_i}(\psi_{h_s}+\sum_{\bigcup_{i=1}^{s-1}h_i\subset h\subset\bar{h}_s}D_{h\bar{h}})^{s-2},  $$
where $h_i$ are such that $h_i\bigcap h_j=\emptyset$, $\forall i\not=j$, and $a, a_i\geq 0$. (see \cite{noi1}, Lemma 4.4)
(They are special because, for a nonzero class $\prod_{i=1}^sD_{h_i\bar{h}_i}$, they are minimal in the $\QQ$-submodule of relations on $\prod_{i=1}^sD_{h_i\bar{h}_i}$ which do not come from a smaller substratum.)

 When $d=3$, the special relation is $\sigma_3\psi_h$ for any $h\in P_2$, hence relation $\tau\sigma_3$. Another important observation comes from comparing equation (2.1) to relation (1) in Lemma 3.1: on $\ol{M}_{h\bar{h}}$, $\psi_h^2=(\psi')^2=\rho$. (this is true for any $d$ odd when $|h|=(d-1)/2$). Hence relation $\rho\sigma_3$, and also, in conjunction with equation (3.1), and Lemma 3.1--(3), $\tau^2=\rho\sigma_2$. Upon reflection, these are all the relations ensuing from Lemma 3.1 and formula (3.1).

\end{proof}

\end{example}


\providecommand{\bysame}{\leavevmode\hbox
to3em{\hrulefill}\thinspace}

\end{document}